\documentclass{article}
\usepackage{spconf,amsmath,graphicx}
\usepackage{amsfonts,amssymb,amsthm}
\usepackage{mathrsfs}
\usepackage{algorithm}
\usepackage{algorithmic}
\usepackage{xcolor}
\usepackage{enumitem}
\usepackage{xspace}
\usepackage{multirow}
\usepackage{hhline}
\usepackage{dsfont}
\usepackage[utf8]{inputenc}

\def\bw{{\mathbf w}}
\def\bx{{\mathbf x}}
\def\by{{\mathbf y}}

\def\bB{{\mathbf B}}
\def\bC{{\mathbf C}}
\def\bD{{\mathbf D}}

\def\bI{{\mathbf I}}

\def\bR{{\mathbf R}}
\def\bW{{\mathbf W}}

\def\ccS{\mathscr{S}}
\def\cD{\mathcal{D}}

\def\RR{\mathbb{R}}
\def\CC{\mathbb{C}}
\def\EE{\mathbb{E}}

\def\CN{\mathscr{C\!N}}
\def\1{\mathds{1}}

\def\bepsilon{\boldsymbol{\epsilon}}
\def\btheta{\boldsymbol{\theta}}

\def\bmu{\boldsymbol{\mu}}

\def\btau{\boldsymbol{\tau}}
\def\bGamma{\boldsymbol{\Gamma}}
\def\bPsi{\boldsymbol{\Psi}}

\def\bzero{\boldsymbol{0}}

\def\defeq{\stackrel{\Delta}{=}}
\def\Fs{F_{\mathsf{s}}}
\def\SNR{{\small\textsf{SNR}}}

\newcommand{\ie}{i.e.\@\xspace}

\renewcommand{\Re}{\operatorname{Re}}

\newtheorem{proposition}{Proposition}
\newtheorem{remark}{Remark}

\title{TIME-SCALE SYNTHESIS FOR LOCALLY STATIONARY SIGNALS}
\twoauthors
  {Adrien Meynard\sthanks{The first author performed this work while at I2M, Aix-Marseille Universit\'e, France.}}
	{Duke University,\\Department of Mathematics,\\
	Durham, NC, USA}
  {Bruno Torr\'esani}
	{Aix Marseille Univ, CNRS, Centrale Marseille,\\
	I2M,\\
	Marseille, France}
\begin{document}
%
\maketitle
\begin{abstract}
We develop a time-scale synthesis-based probabilistic approach for the modeling of locally stationary signals. Inspired by our previous work, the model involves zero-mean, complex Gaussian wavelet coefficients, whose distribution varies as a function of time by time dependent translations on the scale axis. In a maximum \textit{a posteriori} approach, we propose an estimator for the model parameters, namely the time-varying scale translation and an underlying power spectrum. The proposed approach is illustrated on a denoising example. It is also shown that the model can handle locally stationary signals with fast frequency variations, and provide in this case very sharp time-scale representations more concentrated than synchrosqueezed or reassigned wavelet transform.
\end{abstract}
\section{Introduction}
Classical time-frequency analysis is generally used for building signal representations from which relevant information can be extracted (see e.g.~\cite{Carmona98practical,Flandrin18explorations,Grochenig01foundations} for reviews). Under suitable assumptions, linear transforms such as the STFT, wavelet transform or generalizations are invertible, which also leads to so-called synthesis approaches~\cite{Balazs13adapted}. The latter express signals as linear combinations of \textit{time-frequency atoms}, and the corresponding time-frequency coefficients provide another type of time-frequency representation, which is less constrained by consistency requirements and uncertainty principles.

Statistical approaches to time-frequency analysis often rely on \textit{ad hoc} statistical models for time-frequency transforms. Information extraction is then formulated as a statistical estimation problem. Examples include non-negative matrix factorization methods (see~\cite{Fevotte99nonnegative}), detection of time-frequency components~\cite{Reinhold18objective,Huillery08description}, and several other tasks. In most situations, modeling appears as a post-processing stage after computation of a time-frequency transform. However, statistical models are generally not compatible with consistency conditions satisfied by time-frequency transforms.

Recently, Turner and Sahani~\cite{Turner14time} introduced a new Bayesian paradigm, under the name of probabilistic time-frequency representations. The idea is to express signals as the result of a synthesis from a random time-frequency representation, for which a prior distribution is chosen. This is applied to various contexts, such as the synthesis of stationary signals, and several non-stationary situations, including NMF-based component estimation, non-stationary noise. A similar point of view has already been taken by various authors in the past, see for example~\cite{Davy06bayesian} and references therein. We rely here on the Turner-Sahani model, which we revisit in a slightly different way, assuming a generative model of the form
\begin{equation}
\label{fo:obs.eq.cont}
y(t) = \Re\left(\sum_s (\psi_s * W_s)(t)\right)+\epsilon(t)\ ,
\end{equation}
where $\{\psi_s,s\}$ is a filter bank, labeled by a frequency (or scale) index $s$, the $W_s(t)$ are random subband (time-frequency or time-scale) coefficients, and $ \epsilon(t)$ is a noise.

Our focus is here on non-stationary signals, more precisely \textit{locally stationary signals} for which time-dependent spectral characteristics can be defined. Following our earlier JEFAS approach~\cite{Meynard18spectral,Meynard19separation,Omer17time}, we investigate a class of non-stationarity generated by time-dependent shifts in the time-scale domain. JEFAS is analysis based, \ie post-processing of time-scale representation. We introduce JEFAS-S, a Bayesian \textit{synthesis approach}  that provides adaptive time-scale signal representation, together with corresponding parameter estimation. While JEFAS based estimation was based on approximations of the time-scale transform, an exact estimation is possible here, and we provide a corresponding EM algorithm. In addition, in some situations, the model is flexible enough to provide extremely concentrated time-scale representations that can be sharper than reassigned transforms~\cite{Auger95improving}.

\section{The synthesis model}
In this paper, we limit to time-scale representation, \ie wavelet synthesis. We denote by $\psi$ the analysis wavelet, and by $\psi_s$ scaled wavelets defined by $\psi_s(t)=q^{-s/2}\psi\left(q^{-s}t\right)$, for some constant $q>1$. $s \in\{s_1,...\,s_M \}$ is a finite set of scales.

\subsection{The discrete model}
We consider the finite periodic case: assume we have $N$ time values $\btau=(\tau_1\cdots\tau_N)^T$ and the corresponding sampled signal $\by=(y(\tau_1)\cdots y(\tau_N))^T$ with sampling frequency $\Fs$. We then focus on a corresponding discretized wavelet transform. For $n\in\{1,...\, N\}$, $m\in\{1,...\, M\}$, denote by $\psi_{nm}\in\CC^N$ the vector $(\psi_{s_m}(\tau_{1-n})\cdots\psi_{s_m}(\tau_{N-n}))^T$, and by $\bPsi_n\in\CC^{N\times M}$ the matrix obtained by concatenation of vectors $\psi_{nm},m\in\{1,...\, M\}$.
The observation equation~\eqref{fo:obs.eq.cont} then reads 
\begin{equation}
\label{fo:discrete.model}
\by = \by_0+\bepsilon = \Re\left( \sum_{n=1}^{N} \bPsi_n\bw_n\right) + \bepsilon\ ,
\end{equation}
where the $\bw_n\in\CC^M,\,n=1,\dots N$ are vectors of synthesis coefficients. This model can also be written in matrix form as $\by=\Re(\bD\bW +\bepsilon)$, where the dictionary matrix $\bD$ is the concatenation of matrices $\bPsi_n$, and $\bW=\mathsf{vec}(\bw_1,\dots\bw_N)$.

In this paper, $\bepsilon$ will be modeled as a Gaussian white noise, with variance $\sigma^2$, as in~\cite{Turner14time}. Non-stationarity will be introduced via a suitable prior on $\bW$, that intends to describe locally time-warped situations as introduced in~\cite{Clerc03estimating}.

\subsection{A class of non-stationary priors: time warping}
When all subband signals $W_s$ in~\eqref{fo:obs.eq.cont} are stationary, the resulting signal $y$ is stationary. We are interested here in a specific situation where non-stationarity induces a time-dependent shift on the scale axis, as studied in~\cite{Clerc03estimating,Omer17time,Meynard18spectral}. It was shown there that such a model can account for signals obtained by time warping stationary signals, namely signals of the form
\begin{equation}
\label{fo:time.warping}
y(t) = (\cD_\gamma x)(t) \defeq \sqrt{\gamma'(t)}x(\gamma(t))\ ,
\end{equation}
where $x$ is a wide sense stationary random signal, and $\gamma$ is a smooth, strictly increasing function.

To build the prior distribution on discrete subband coefficients, we make the following assumptions
\begin{itemize}
\vspace{-2mm}
\item
The vectors $\bw_n$ are decorrelated, zero-mean, circular complex Gaussian vectors: $\bw_n\sim\CN_c(\mathbf{0},\bC_n)$
\vspace{-3mm}
\item
The corresponding covariance matrices $\bC_n$ are translates of a fixed function $f$ as shown in~\cite{Meynard18spectral}, namely
\begin{equation}
\label{fo:cov.function}
[\bC_{n}]_{mm'} \defeq [\bC(\theta_n)]_{mm'} = f(s_m\!+\!\theta_n,s_{m'}\!+\!\theta_n)\ ,
\end{equation}
where $f:\left(\RR_+^*\right)^2\rightarrow\CC$ is Hermitian and positive-semidefinite, and $\theta_n\in\RR$ is the shift parameter.
\end{itemize}
In~\cite{Meynard18spectral}, it was shown that the wavelet coefficients of a stationary random signal modified by time warping transform can be approximated by random vectors satisfying the above assumptions. There, the parameter $\theta_n$ represents a local dilation factor at time $\tau_n$ (derivative $\gamma'(\tau_n)$ of the time warping function at $\tau_n$), and $f$ involves the power spectrum $\ccS$ of the underlying signal and the Fourier transform of the wavelet:
\begin{equation}
\label{fo:cov.function.wav}
f(s,s') = q^{\frac{s+s'}{2}}\int_0^\infty \ccS(\xi)\overline{\hat\psi(q^s\xi)}\hat\psi(q^{s'}\xi)\,d\xi\ .
\end{equation}

\section{Estimation procedure}
\subsection{Bayesian inference}
The estimation of the subband coefficient matrix $\bW$ relies on the evaluation of the corresponding posterior distribution. The latter depends on the following parameters, which are supposed to be known at this point: the dilation factors $\theta_n$, and the covariance function $f$.
Let $\bGamma_0\in\CC^{MN\times MN}$ be the block diagonal matrix with blocks $\bC_1,\dots\bC_N$. The posterior distribution of the subband coefficient is a complex Gaussian law $p(\bW|\by)\sim \CN(\bmu,\bGamma,\bR)$, with mean and covariance
\[
\bmu = \bGamma_0\bD^H \bC_y^{-1}\by\ ,\qquad
\bGamma = \bGamma_0 - \dfrac1{4}\bGamma_0\bD^H\bC_y^{-1}\bD\bGamma_0\ ,
\]
(the relation matrix $\bR$, not useful here is not provided) where
\begin{equation}
\bC_y = \sigma^2\bI + \dfrac12\Re\left(\bD\bGamma_0\bD^H\right)\ .
\end{equation}
Therefore the posterior expectation $\tilde\bw_n$ of $\bw_n$ reads
\begin{equation}
\label{fo:transform.adap}
\tilde\bw_n  = \dfrac12\bC_n\bPsi_n^H\bC_y^{-1}\by\ ,
\end{equation}
where the matrix $\bC_y$ can be expressed as
\begin{equation}
\label{fo:transform.adap.2}
\bC_y = \sigma^2\bI + \dfrac12\Re\left(\sum_{n=1}^{N}\bPsi_n\bC_n\bPsi_n^H\right)\ .
\end{equation}
\begin{remark}
It is worth mentioning that unlike the prior distribution, the posterior distribution of subband coefficients involves time correlations. Indeed, given any $n,n'=1\dots N$,
\begin{equation*}
\EE\{\tilde\bw_n\tilde\bw_{n'}^H\} = \delta_{nn'}\bC_n -\dfrac1{4}\bC_n\bPsi_n^H\bC_y^{-1}\bPsi_{n'}\bC_{n'}\ ,
\end{equation*}
which generally does not vanish when $n'\neq n$.
\end{remark}

\subsection{Parameter selection and reconstruction}
We now discuss the choice of the two model parameters, \ie the scaling factors $\theta_n$ and the covariance function $f$. We first notice that the expression in equation~\eqref{fo:cov.function.wav} provides a natural choice for the covariance function $f$. The latter involves the wavelet, which is known, and the power spectrum $\ccS$ of the underlying stationary process, which is unknown. In this setting, we then have to provide the vector $\btheta=(\theta_1,\dots\theta_N)$ of scaling factors and the power spectrum.

The JEFAS algorithm introduced in~\cite{Meynard18spectral} provides a valuable, analysis-based approach for estimating $\btheta$ and $\ccS$. We now describe an alternative algorithm, based on the EM (Expectation Maximization) principle~\cite{Dempster77maximum}, especially tailored for the synthesis approach developed in this paper. Here, $\by$ is the observation, $\btheta$ the parameter, and $\bW$ the latent variable.

\begin{proposition}[EM steps]
\label{prop:EM}
Denote by $\tilde\btheta^{(k-1)}$ the vector of dilation factors at iteration $k-1$ of the algorithm. Let $\tilde\bW^{(k-1)}$ be the matrix of subband coefficients at iteration $k-1$. Then the update at iteration $k$ relies on the following two steps:
\begin{enumerate}
\item For $n\in\{1,\ldots,N\}$, the estimate~\eqref{fo:transform.adap} at time $\tau_n$ reads

\vspace{-5mm}
\begin{equation}
\tilde\bw_n^{(k)}  = \dfrac12\bC\left(\tilde\theta_n^{(k-1)}\right)\bPsi_n^H\bC_y\left(\tilde\btheta^{(k-1)}\right)^{-1}\by\ .
\label{fo:wav.EM}
\end{equation}

\vspace{-5mm}
\item The scaling factor $\tilde\btheta$ is re-estimated by solving
\begin{equation}
\label{fo:pb.EM}
\tilde\theta_n^{(k)}=\arg\min_{\theta}\ Q_{kn}(\theta)\ ,
\end{equation}
\vspace{-6mm}
\begin{eqnarray*}
Q_{kn}(\theta)\!\!\!&\!\!=\!\!\Big[\log\! |\!\det(\bC(\theta))|\! +\! \tilde\bw_n^{(k)H}\bC(\theta)^{-1}\tilde\bw_n^{(k)}\\
&+\mathrm{Trace}\left(\bC(\theta)^{-1}\bGamma_{n}\!\!\left(\tilde\btheta^{(k-1)}\right)\!\right)\!\Big]\ ,
\end{eqnarray*}
$\bGamma_{n}\left(\tilde\btheta^{(k-1)}\right)\in\CC^{M\times M}$ being the $n$-th diagonal block of the posterior covariance matrix $\bGamma\left(\tilde\btheta^{(k-1)}\right)$. 
\end{enumerate}
\end{proposition}
After running the corresponding algorithm (described in more details below), an estimate for the time-scale coefficients $\tilde\bw$ is available, and a corresponding estimate $\tilde\by_0$ for the signal $\by_0$ can be obtained as

\vspace{-3mm}
\begin{equation}
\label{fo:reco}
\tilde\by_0 = \Re\left( \sum_{n=1}^{N} \bPsi_n\tilde\bw_n\right)\ .
\end{equation}
\vspace{-3mm}

Notice that the reconstruction expression~\eqref{fo:reco} combined with~\eqref{fo:transform.adap} can be interpreted as a Wiener filtering. The bias and variance of the estimator can be evaluated.
\begin{proposition}
With the above notation, the bias of the estimator $\tilde\by_0$ is given by
\begin{equation}
\bB \defeq\EE\left\lbrace \tilde\by_0|\by_0 \right\rbrace - \by_0 = -\sigma^2\bC_y^{-1}\by_0\ ,
\end{equation}
and the corresponding error variance reads
\begin{align}
\nonumber
\bR(\tilde\by_0|\by_0) &\defeq \EE\left\lbrace \left. \left(\tilde\by_0 - \EE\{\tilde\by_0|\by_0\}\right) \left(\tilde\by_0 - \EE\{\tilde\by_0|\by_0\}\right)^T \right| \by_0 \right\rbrace \\
&= \sigma^2\left(\bI-\sigma^2\bC_y^{-1}\right)^2\ .
\end{align}
\end{proposition}

\subsection{Algorithm: JEFAS-Synthesis}

The steps of the estimation algorithm are given in Algorithm~\ref{alg:JEFAS.adap}. The latter takes as input the signal $\by$, the noise variance $\sigma^2$, a precision parameter $\Lambda$ for the stopping criterion and a bandwidth parameter $N'$ (see below).

\textit{Initialization.} The algorithm requires initial estimates $\btheta^{(0)}$ for the parameters, and the function $f$ in~\eqref{fo:cov.function}. In JEFAS-S, we use the expression~\eqref{fo:cov.function.wav}, for which an initial estimate of $\ccS$ has to be provided. When successful, JEFAS~\cite{Meynard18spectral} provides such an estimate. Otherwise, a rough estimate can be obtained from the Welch periodogram of the input signal $\by$.

\textit{Stopping criterion.} EM guarantees the monotonicity of the Likelihood function $\mathcal{L}(\btheta)$. The increment of the latter is used as a stopping criterion: EM will stop when the condition
\begin{equation}
\label{fo:crit.jefass}
\mathcal{L}(\btheta^{(k)})-\mathcal{L}(\btheta^{(k-1)})<\Lambda
\end{equation}
is true. Here $\Lambda>0$ is a parameter fixed by the user.

\textit{Dimension reduction.}
The matrix $\bC_y$ of dimension $NM\times NM$ can be extremely large. However, it generally has fast off-diagonal decay. This can be exploited to speed up the evaluation of $\tilde \bw_n$ in~\eqref{fo:wav.EM} by restricting to a neighborhood $[n-N'/2,n+N'/2]$ of $n$ of given bandwidth $N'$.

\textit{Optimization.}
The optimization problem~\eqref{fo:pb.EM} is solved using a standard quasi-Newton scheme.

\textit{Spectrum estimate update.}
The spectrum update from the current estimate of $\bW$ is performed in two steps: first correct for the translation by $\theta_n$, to obtain an approximately stationary subband transform, then average over time to obtain a wavelet based spectral estimate as in~\cite{Meynard18spectral}.

\begin{algorithm}[t]
\caption{$(\tilde\bW,\tilde\btheta,\tilde\ccS_X)=\operatorname{JEFAS-S}(\by,\sigma^2,\Lambda,N')$}
\begin{algorithmic}
\STATE $\bullet$ {\bfseries Initialization:} estimate $\tilde\btheta^{(0)}$ and $\tilde\ccS^{(0)}$ using JEFAS.
\STATE $\bullet$ $k\leftarrow 1$.
\WHILE{stopping criterion~\eqref{fo:crit.jefass} = \textsf{FALSE}}
\FORALL{$n\in\{1,\ldots,N\}$}
\STATE $\bullet$ Restrict $\bPsi_n$, $\bC_y\left(\tilde\btheta^{(k-1)}\right)$ and $\by$ to the interval  $[n-N'/2,n+N'/2]$.
\STATE $\bullet$ Compute $\tilde\bw_n^{(k)}$ using~\eqref{fo:wav.EM}.
\ENDFOR
\STATE $\bullet$ Estimate $\tilde\btheta^{(k)}$ by solving~\eqref{fo:pb.EM}.
\STATE $\bullet$ Estimate $\tilde\ccS^{(k)}$ using the wavelet based estimate.
\STATE $\bullet$ $k\leftarrow k+1$.
\ENDWHILE
\end{algorithmic}
\label{alg:JEFAS.adap}
\end{algorithm}

\begin{remark}
Other choices can be made for the function $f$, which can lead to different estimates for subband coefficients, while preserving reconstruction (see section~\ref{sse:JEFASS.sig.semi.real}).
\end{remark}

\section{Numerical results}

\subsection{Illustration on denoising of a synthetic signal}
\label{sse:JEFASS.toy.sig}

We first evaluate the performances of JEFAS-S on a denoising problem. A synthetic non-stationary signal $\by$ is built as follows: start from a stationary signal $\bx$, with power spectrum $\ccS$ equal to the sum of two non-overlapping Hann windows, and apply the time warping deformation $\cD_\gamma$ to $\bx$, with $\gamma'$ an exponentially damped sine wave. Here, $\bx$ is one second long, sampled at $\Fs=8192$~Hz. 

We denote by \SNR$_\by$ and \SNR$_{\tilde\by_0}$ the input and output signal-to-noise ratios. Numerical results show that \SNR$_{\tilde\by_0}$ is larger than \SNR$_\by$ as long as \SNR$_\by$ is in the range $[2~\hbox{dB}, 25~\hbox{dB}]$, with maximal improvement of $8$~dB.
The $25$~dB upper limit for \SNR$_{\tilde\by_0}$ is presumably due the distortion intrinsically introduced by the reconstruction formula~\eqref{fo:reco}: bias and variability in the time warping estimation.

\begin{figure}
\centering
\includegraphics[width=.47\textwidth]{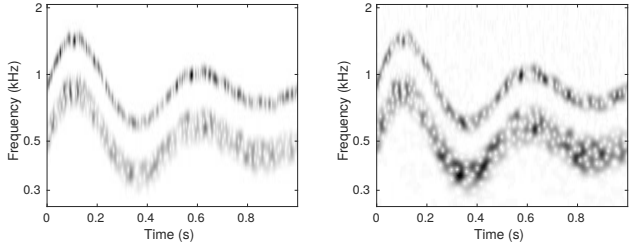}
\vspace{-3mm}
\caption{Synthetic signal. Left: representation given by JEFAS-S. Right: scalogram (wavelet transform).}
\label{fig:JEFASS.temps.echelle}
\vspace{-2mm}
\end{figure}

In the specific case where the input \SNR\xspace is $16$~dB, and after initializing with the output of JEFAS, JEFAS-S converges in $3$ iterations (CPU time: 347 seconds on a computer running an Intel Xeon E5-2680 v4 processor). JEFAS-S does not significantly improve the quality of the estimated time warping function. Indeed, the mean square error on the time warping function estimation decreases by about $0.5\%$ from JEFAS to JEFAS-S. We display the estimated adapted time-scale representation $\tilde\bW$ in Fig.~\ref{fig:JEFASS.temps.echelle} (left). As expected, it is very similar to the wavelet transform (right), though a bit sharper. Indeed, the choice of the expression~\eqref{fo:cov.function.wav} for the covariance function $f$ yields a wavelet-like representation. The main visible difference concerns the temporal oscillations of $\tilde\bW$, due to the prior assumption of temporal decorrelation between $\bw_n$.

\subsection{Locally harmonic signal with fast varying frequency}
\label{sse:JEFASS.sig.semi.real}
\vspace{-1mm}
We now consider a locally harmonic signal, of the form
\[
y(t) = A(t)\cos(2\pi\phi(t))\ ,
\vspace{-1mm}
\]
where the instantaneous frequency $\phi'$ is a fast varying function chosen as the measurement of the heart rate of a person suffering from atrial fibrillation (real data). The synthetic instantaneous amplitude $A$ is a slowly varying function, the signal is termed ``semi-real''. Apart from the amplitude modulation, this signal follows the model~\eqref{fo:time.warping}: the time warping function derivative is the instantaneous frequency, and the underlying stationary signal $x$ is sinusoidal. The signal duration is 83.1~seconds, sampled at $\Fs=10$~Hz ($N=832$ samples).
 
\begin{figure}
\centering
\includegraphics[width=.42\textwidth]{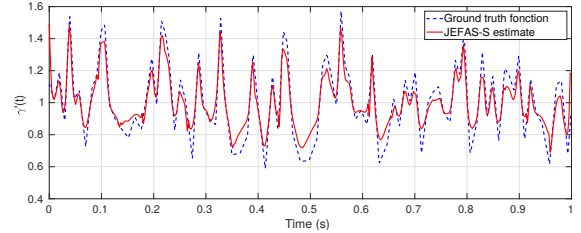}
\vspace{-4mm}
\caption{Semi-real signal. Estimated time warping function compared with the normalized instantaneous frequency $\phi'$.}
\label{fig:JEFASS.tw.sigreel}
\end{figure}

Because of the fast instantaneous frequency variations, the wavelet transform of $y$ (not shown here) contains interference patterns, the model in~\cite{Meynard18spectral} is not adequate and JEFAS does not converge. We initialized JEFAS-S to $\btheta=\bzero$, and a constant function for $\ccS$. Given these initial values (far from actual values), JEFAS-S converges slowly (72 iterations). Results are displayed in Fig.~\ref{fig:JEFASS.tw.sigreel}, where the estimated instantaneous frequency is superimposed on the ground-truth. This shows that JEFAS-S is indeed able to estimate fast varying frequency modulations.

\begin{figure}
\centering
\includegraphics[width=.42\textwidth]{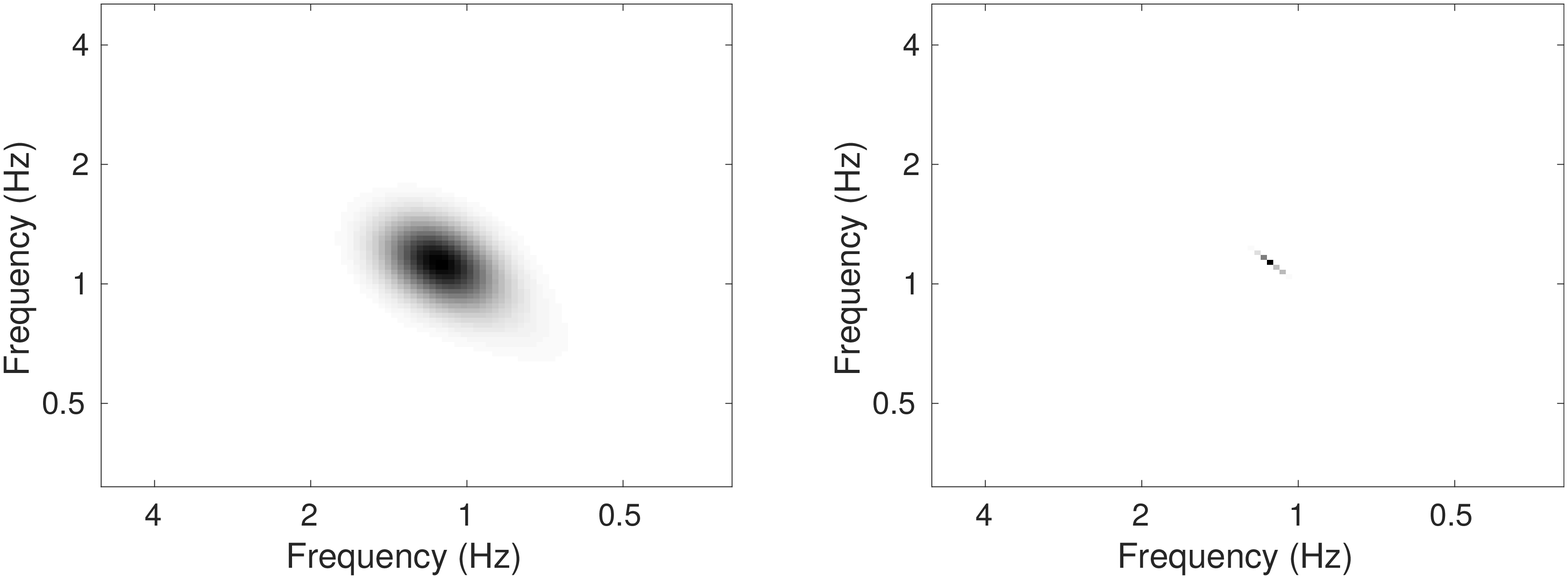}
\includegraphics[width=.42\textwidth]{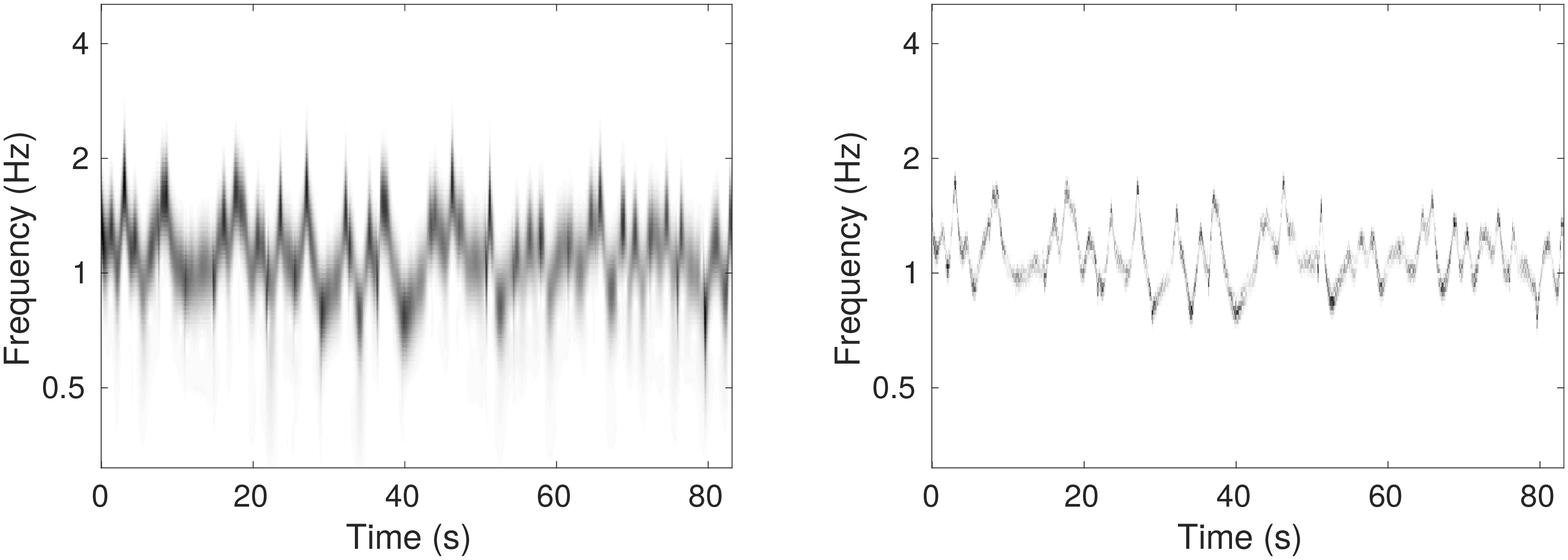}
\vspace{-4mm}
\caption{Semi-real signal. Top: two priors for the covariance matrix. Bottom: associated time-scale representations.}
\label{fig:comp.f.function}
\end{figure}

In addition, one can take advantage of this result to obtain a sharper time-scale representation. To that end, we choose a sharply concentrated prior covariance function $f_\sharp$, of the form $f_\sharp(s,s') = \exp\{-(s-\varsigma)^2/\sigma_s^2\}\delta_{ss'}$, where $\nu_1$ denotes the central frequency of the sine wave, $\varsigma = \log_q\left(\xi_0/\nu_1\right)$, and $\sigma_s$ is a tuning parameter for the scale concentration. We display in the top of Fig.~\ref{fig:comp.f.function} the covariance matrices $\bC(0)$ corresponding to the expression~\eqref{fo:cov.function.wav} (left) and $f_\sharp$ (right), which is indeed very sharp. The corresponding estimated time-scale representations are displayed on the the bottom images of Fig.~\ref{fig:comp.f.function}.
The new prior is clearly adapted to locally harmonic signals, \ie signals with a sparse underlying spectrum. Thus, in such situations, JEFAS-S enables the construction of sharp time-scale representations, competing with standards techniques such as synchrosqueezing. Furthermore, we stress that the quality of the reconstruction is not degraded.

\vspace{-1mm}
\section{Conclusion}
\vspace{-2mm}
We have described an alternative to the JEFAS model of~\cite{Meynard18spectral} for locally deformed signals. Unlike JEFAS, which is an analysis based approach (\ie post-processing of wavelet transform), JEFAS-S is synthesis-based and therefore less constrained by uncertainty principles. We illustrated the JEFAS-S on a denoising example. Our numerical results also show that JEFAS-S is able to handle locally stationary signals with fast varying instantaneous frequency, and can provide very sharp time-scale representations.

%
%

While the current paper was focused on wavelet transform, the JEFAS-S model can handle arbitrary subband decompositions (such as the NSDGT~\cite{Balazs11theory}). Such extensions will be discussed in a forthcoming publication, together with additional numerical results and complete proofs. JEFAS-S can also be extended to more general transformations, for example involving amplitude modulations or filtering posterior to time warping. This is an ongoing work.

\bibliographystyle{IEEEbib}
\bibliography{ICASSP20}

\begin{thebibliography}{10}

\bibitem{Carmona98practical}
Ren\'e Carmona, Wen-Liang Hwang, and Bruno Torr\'esani,
\newblock {\em Practical time-frequency analysis: Gabor and Wavelet Transforms
  With an Implementation in S},
\newblock Academic Press, 1998.

\bibitem{Flandrin18explorations}
Patrick Flandrin,
\newblock {\em Explorations in Time-Frequency Analysis},
\newblock Cambridge University Press, 2018.

\bibitem{Grochenig01foundations}
Karlheinz Gr{\"o}chenig,
\newblock {\em Foundations of time-frequency analysis},
\newblock Applied and Numerical Harmonic Analysis. Birkh\"auser Inc., Boston,
  MA, 2001.

\bibitem{Balazs13adapted}
Peter Balazs, Monika D{\"o}rfler, Matthieu Kowalski, and Bruno Torr{\'e}sani,
\newblock ``{Adapted and adaptive linear time-frequency representations: a
  synthesis point of view},''
\newblock {\em {IEEE Signal Processing Magazine}}, vol. 30, no. 6, pp. 20--31,
  Nov. 2013.

\bibitem{Fevotte99nonnegative}
C\'edric F\'evotte and Ali~T. Cemgil,
\newblock ``Nonnegative matrix factorisations as probabilistic inference in
  composite models,''
\newblock in {\em Proc.~17th European Signal Processing Conference (EUSIPCO)},
  Glasgow, Scotland, Aug. 2009, pp. 1913--1917.

\bibitem{Reinhold18objective}
Isabella Reinhold, Maria Sandsten, and Josefin Starkhammar,
\newblock ``Objective detection and time-frequency localization of components
  within transient signals,''
\newblock {\em The Journal of the Acoustical Society of America}, vol. 143, no.
  4, pp. 2368--2378, 2018.

\bibitem{Huillery08description}
Julien Huillery, Fabien Millioz, and Nadine Martin,
\newblock ``On the description of spectrogram probabilities with a chi-squared
  law,''
\newblock {\em IEEE Transactions on Signal Processing}, vol. 56, no. 6, pp.
  2249 -- 2258, June 2008.

\bibitem{Turner14time}
Richard~E. {Turner} and Maneesh {Sahani},
\newblock ``Time-frequency analysis as probabilistic inference,''
\newblock {\em IEEE Transactions on Signal Processing}, vol. 62, no. 23, pp.
  6171--6183, Dec. 2014.

\bibitem{Davy06bayesian}
Manuel Davy, Simon~J. Godsill, and J{\'e}r{\^o}me Idier,
\newblock ``{Bayesian Analysis of Polyphonic Western Tonal Music},''
\newblock {\em {Journal of the Acoustical Society of America}}, vol. 119, no.
  4, pp. 2498--2517, 2006.

\bibitem{Meynard18spectral}
Adrien Meynard and Bruno Torr{\'e}sani,
\newblock ``{Spectral Analysis for Nonstationary Audio},''
\newblock {\em {IEEE/ACM Transactions on Audio, Speech and Language
  Processing}}, vol. 26, no. 12, pp. 2371 -- 2380, Dec. 2018.

\bibitem{Meynard19separation}
Adrien Meynard,
\newblock ``{S{\'e}paration de sources doublement non stationnaire},''
\newblock in {\em {GRETSI 2019 - XXVII{\`e}me Colloque francophone de
  traitement du signal et des images}}, Lille, France, Aug. 2019.

\bibitem{Omer17time}
Harold Omer and Bruno Torr{\'e}sani,
\newblock ``{Time-frequency and time-scale analysis of deformed stationary
  processes, with application to non-stationary sound modeling},''
\newblock {\em {Applied and Computational Harmonic Analysis}}, vol. 43, no. 1,
  pp. 1 -- 22, 2017.

\bibitem{Auger95improving}
François {Auger} and Patrick {Flandrin},
\newblock ``Improving the readability of time-frequency and time-scale
  representations by the reassignment method,''
\newblock {\em IEEE Transactions on Signal Processing}, vol. 43, no. 5, pp.
  1068--1089, May 1995.

\bibitem{Clerc03estimating}
Maureen Clerc and Stéphane Mallat,
\newblock ``Estimating deformations of stationary processes,''
\newblock {\em Ann. Statist.}, vol. 31, no. 6, pp. 1772--1821, Dec. 2003.

\bibitem{Dempster77maximum}
Arthur~P. Dempster, Nan~M. Laird, and Donald~B. Rubin,
\newblock ``{M}aximum {L}ikelihood from {I}ncomplete {D}ata via the {EM}
  {A}lgorithm,''
\newblock {\em Journal of the Royal Statistical Society. Series B
  (Methodological)}, vol. 39, no. 1, pp. 1--38, 1977.

\bibitem{Balazs11theory}
Peter Balazs, Monika Dörfler, Florent Jaillet, Nicki Holighaus, and Gino
  Velasco,
\newblock ``Theory, implementation and applications of nonstationary {G}abor
  frames,''
\newblock {\em Journal of Computational and Applied Mathematics}, vol. 236, no.
  6, pp. 1481 -- 1496, 2011.

\end{thebibliography}

\end{document}